\documentclass[conference]{IEEEtran}
\usepackage{amsmath}
\usepackage{pdflscape}
\usepackage[justification=centering]{caption}
\usepackage{ amssymb }
\newtheorem{thm}{Theorem}[section]
\newtheorem{prop}[thm]{Proposition}

\newtheorem{conj}{Conjecture}
\newcommand{\pf}{{\bf Proof. \ }}
\newcommand{\qed}{\hfill $\Box$ \\}
\renewcommand{\baselinestretch}{1.15}

\begin{document}

\title{The joint weight enumerator of an LCD code and its dual }
\author
{Adel Alahmadi$^{1}$, Michel Deza$^{2}$, Mathieu Dutour-Sikiri\'c$^{3}$, Patrick Sol\'e$^{1,4}$\\
\small
$^1$ Math Dept.,King Abdulaziz University, Jeddah, Saudi Arabia\\ \small
$^2$ University of Campinas, IMECC, Campinas, Brazil\\ \small
$^3$ Rudjer Boskovi\'c Institute, Bijenicka 54, 10000, Zagreb, Croatia\\ \small
$ ^4$LTCI, Telecom ParisTech, Paris, France
}
\date{}
\maketitle
\begin{abstract} 
A binary linear code is called {\em LCD} if it intersects its dual trivially.
We show that the coefficients of the joint weight enumerator of such a code with its dual
satisfy linear constraints, leading to a new linear programming bound on the size of an LCD code of given length and minimum distance. In addition, we show that
this polynomial is, in general, an invariant of a matrix group of dimension $4$ and order $12$.
Also, we sketch a Gleason formula for this weight enumerator.
\end{abstract}

\section{Introduction} 
A binary linear code is called {\em LCD} if it intersects its dual trivially.
These codes, introduced by Massey in \cite{Massey}, give an
optimum linear coding solution for the two user binary adder
channel.
They were rediscovered recently  in \cite{CarGui}
in a context of countermeasures to passive and active side channel analyses on
embedded crytosystems. While most studies so far are concerned with constructions,
the recent article \cite{DKOSS} contains a linear programming bound on the size of
a binary linear LCD code of given length and minimum distance.
This bound is proved there to be sharper in many instances than the classical
linear programming bound.

In the present work, we will present a linear programming bound on the same quantity
with variables being the coefficients of the joint weight enumerator of such a code with its dual.
The advantage of using this four-variable polynomial is that the condition of LCD-ness
is now linear (instead of quadratic in \cite{DKOSS}) and necessary and sufficient
(instead of only necessary in \cite{DKOSS}).

In addition, we observe that the joint weight enumerator of a linear code and its dual
is an invariant of a matrix group of dimension $4$ and order $12$,
and sketch a Gleason formula for this weight enumerator.
This seems to have been unnoticed since the seminal application of invariant theory to
codes $50$ years ago in \cite{MMS}.

The material is organized as follows. Section II collects the notation and the definitions
needed in the following sections.
Section III, of independent interest develops the invariant theory of the joint weight
enumerator of a code and its dual. Section III establishes the linear programming bound
with variables the coefficients of the joint weight enumerator of a LCD code and its dual,
and validates it by improving the numerical results of \cite{DKOSS}.

\section{Definitions and notation}

\subsection{Codes}
In this work, all considered codes are binary and linear.
A {\em code of length} $n$ is thus subspace over the field $GF(2)$ of the vector space $GF(2)^n$.
The {\em dual} $C^\bot$ of such a code $C\le GF(2)^n$ is defined 
with respect of the standard euclidean inner product as
$$C^\bot=\{ y \in GF(2)^n |\, \forall x\in C,\, x.y=0\}.$$
A binary linear code $C$ is {\em Linear Complementary Dual} (LCD) if $C\bigcap C^\bot=0$.

\subsection{Weight distributions}
The {\em Hamming weight} of a binary vector $x$ of length $n$ is the number
of its nonzero coordinates.
For $i=0,1,\cdots,n,$ the {\em weight distribution} $A_i$ of a binary code
of length $n$ is the number of code words with weight $i$.

The weight enumerator of a code $C$ is then
\begin{equation*}
w_C(x,y) = \sum_{i=0}^n A_i x^{n-i} y^{i}
\end{equation*}

\subsection{Joint weight enumerator}
Let $u,v$ denote binary vectors of length $n$.
We define $i(u,v)$, $j(u,v)$, $k(u,v)$ and $l(u,v)$ to be the number of indices
$i\in \{1,n\}$ with $(u_i, v_i) = (0, 0)$, $(0,1)$, $(1,0)$ and $(1,1)$, respectively.

The {\em joint weight enumerator} $J(A,B)$ of, say,  two binary linear codes $A,B$, is the four-variable polynomial
 defined by the formula
\begin{equation}
J(A,B)(a,b,c,d)=\sum_{u\in A,\, v\in B}  a^{i(u,v)} b^{j(u,v)} c^{k(u,v)} d^{l(u,v)} .
\end{equation}
Regrouping the terms on the basis of monomials, we get
\begin{equation}
J(A,B)(a,b,c,d)=\sum_{i+j+k+l = n} M(i,j,k,l) a^{i} b^j c^ k d^l.
\end{equation}
Note that the number of $M(j,k,l)'$s is ${n-1\choose 3}$, the number of {\em decompositions} (ordered partitions)
of the integer $n$ into four parts.

\section{Invariant theory of the joint weight enumerator}\label{INV_section}
In this Section, the code $C$ needs not to be LCD.
{\prop 
\label{formula1}
$$ J(C,C^{\bot})(a,b,c,-d)=J(C,C^{\bot})(a,b,c,d). $$}

\pf
By orthogonality of $u \in C$ and $v \in C^{\bot}$, we see that $l$ is even.
\qed

{\prop 
\label{formula2}
\begin{equation}
\begin{array}{c}
  J(C,C^{\bot})(a,b,c,d)=\frac{1}{2^n} J(C,C^{\bot})(a+b+c+d,\\
                        a+b-c-d,a-b+c-d,a-b-c+d).
\end{array}
\end{equation}
}
\pf
Combine MacWilliams identity \cite[(32) p.148]{TheoryErrorCorrectI} 
between $J(C,C^{\bot})$ and $J(C^{\bot},C)$ with the relation 
\cite[(29) p.148]{TheoryErrorCorrectI}.
\qed
Further, there is the following relation for $n$ even.
{\prop If $n$ is even, then $$ J(C,C^{\bot})(-a,-b,-c,-d)=J(C,C^{\bot})(a,b,c,d). $$}

\pf
Follows by homogeneity of the polynomial.
\qed

The polynomial $J(C,C^{\bot})$ is an invariant of degree $n$ of a group $G=\langle H,J,-I\rangle$ of order $12,$ where 
\begin{equation}
2H=\begin{pmatrix}
1 &  1   & 1   & 1  \\
1 &  1   & -1  & -1 \\
1 &  -1  & 1   & -1 \\
1 & - 1  & -1  & 1
\end{pmatrix},
\end{equation}
and 
\begin{equation}
J=\begin{pmatrix}
1 &  0  &  0  & 0  \\
0 &  1  &  0  & 0  \\
0 &  0  &  1  & 0  \\
0 &  0  &  0  & -1
\end{pmatrix}.
\end{equation}

The three generators  are implied by the  three above
 propositions. The Molien series of $G$ is
\begin{equation}
\frac{1+2t^2+t^4}{(1-t^2)^3(1-t^6)}.
\end{equation}
Denoting by ${\mathcal A}_1$, ${\mathcal A}_2$, ${\mathcal A}_3$, ${\mathcal A}_4$ the primary invariants of respective degrees $2,2,2,6$ and by ${\mathcal B}_1$, ${\mathcal B}_2$, ${\mathcal B}_3$ secondary invariants of respective degree $2,2,4$, we obtain
the following {\em Gleason formula} for even $n\ge 4$:

\begin{equation}\label{ExplicitExpression}
\begin{array}{l}
  J(C,C^{\bot})(a,b,c,d) = \\
\hspace{0.5cm}\sum_{2a_1+2a_2+2a_3+4a_4=n}\alpha_{a_1,a_2,a_3,a_4}\prod_{i=1}^4 {\mathcal A}_i^{a_i}\\
\hspace{0.5cm}+ {\mathcal B}_1 \sum_{a_1+2a_2+2a_3+4a_4=n-2}\beta_{a_1,a_2,a_3,a_4}\prod_{i=1}^4 {\mathcal A}_i^{a_i}\\
\hspace{0.5cm}+ {\mathcal B}_2 \sum_{a_1+2a_2+2a_3+4a_4=n-2}\gamma_{a_1,a_2,a_3,a_4}\prod_{i=1}^4 {\mathcal A}_i^{a_i}\\
\hspace{0.5cm}+ {\mathcal B}_3 \sum_{a_1+2a_2+2a_3+4a_4=n-4}\delta_{a_1,a_2,a_3,a_4}\prod_{i=1}^4 {\mathcal A}_i^{a_i},
\end{array}
\end{equation}
where  $\alpha'$s and $\beta$, $\gamma$, $\delta$'s are arbitrary rational constants.
These invariants are too large to be displayed here.

\section{Linear programming bound}
Let us consider the problem of existence of an LCD code and how linear programming can help.
Given a set of parameters $[n, k, d]$, we will consider a set of equalities and inequalities satisfied by
the parameters $M(i,j,k,l)$.
Our approach is different but related to the one \cite{DKOSS}, where the size of the code is the objective
function and the problem is nonlinear and depends only on $n$ and $d$.
Here, for each triple $[n,k,d]$ we have a convex body defined by it and we stay in the framework of linear
equalities and inequalities. In order to test if this convex body is empty or not, we have to solve a number
of linear programs.

\subsection{Formulation of the problem}\label{EquationSystem}

Let us take a LCD code $C$ and consider the joint LCD enumerator 
\begin{equation}
J(C,C^{\perp})(a,b,c,d) = \sum_{i + j + k + l=n} M(i,j,k,l) a^{i} b^j c^ k d^l.
\end{equation}
The invariance property of Section \ref{INV_section} can be used to find a set of equations
satisfied by the joint LCD enumerator. Alternatively, one can use the basis found 
in Equation \eqref{ExplicitExpression} and compute with it thereafter.
Additionally, we consider the weight distributions $A_i$ and $B_i$ of $C$ and $C^{\perp}$ to be parts
of the equation system.

We have following relations on the coefficients:
\begin{enumerate}
\item Since the dimension of the code $C$ is $k$, we have
\begin{equation}
2^k = \sum_{k=0}^n M(n-k,0,k,0).
\end{equation}
\item Since there is no vector $u\in C$, $v\in C^{\perp}$ with $u.v=1$, we  have
\begin{equation}
0 = M(i,j,k,l) \mbox{~if~} l\equiv 1\pmod 2.
\end{equation}
\item All coefficients are non-negative $M(i,j,k,l) \geq 0$.
\item The weight enumerator constraint is
\begin{equation}
A_i = M(n-i,0,i,0) \mbox{~for~}0\leq i\leq n.
\end{equation}
\item The dual code enumerator constraint is
\begin{equation}
B_i = M(n-i,i,0,0) \mbox{~for~}0\leq i\leq n.
\end{equation}
\item The code packing argument gives
\begin{equation}
A_i + B_i \leq {n\choose i}.
\end{equation}
\item Since $C\cap C^{\perp} = 0$, we  have
\begin{equation}
M(i,0,0,n-i)=0.
\end{equation}
for $0\leq i\leq n-1$.
\item The minimal distance of the code is $d$; so, we have
\begin{equation}
A_i = 0
\end{equation}
for $0 < i < d$.
\item The number of pairs $(u,v)$ with $u$ of weight $p$ is 
\begin{equation}
A_p 2^{n-k} = \sum_{i+j+k+l=n, k+l=p} M(i,j,k,l)
\end{equation}
for $0\leq p\leq n$.
\item The number of pairs $(u,v)$ with $v$ of weight $p$ is 
\begin{equation}
B_p 2^{k} = \sum_{i+j+k+l=n, j+l=p} M(i,j,k,l)
\end{equation}
for $0\leq p\leq n.$
\end{enumerate}
We denote by ${\mathcal K}(n,k,d)$ the polytope determined by all those equalities and inequalities.

\subsection{Relation with the simplified problem of \cite{DKOSS}}
In \cite{DKOSS}, a formulation of linear programming bound is given for the coefficients $A_i$ and $B_i$
given above.
For the sake of completeness, we reformulate it
 in our language of polytopes, since the
formulation of \cite{DKOSS} is very
 different.

The coefficients $A_i$ and $B_i$ satisfy the following constraints:
\begin{enumerate}
\item $A_i\geq 0$ and $B_i\geq 0$ for $1\leq i\leq n$,
\item $A_0=B_0=1$,
\item $A_i=0$ for $1\leq i< d$,
\item for $1\leq i\leq n$ the inequality
\begin{equation}\label{BinomialEq_res}
A_i + B_i \leq {n\choose i},
\end{equation}
\item for $0\leq i\leq n$ we have
\begin{equation}\label{RelationA_Bi}
B_i = 2^{-k} \sum_{j=0}^n A_j P_i(j)
\end{equation}
with $P_i$ being the $i$-th Krawtchuk polynomial defined by the relation
\begin{equation}
\sum_{i=0}^n P_i(x) z^i = (1+z)^{n-x} (1-z)^x.
\end{equation}

\end{enumerate}
Note that the inequality $B_i\geq 0$ is 
 the {\em Delsarte inequality} and 
  the only equation, that is specific to LCD codes, is \eqref{BinomialEq_res}.
All those linear equalities and inequalities define a polytope ${\mathcal K}_{res}(n,k,d)$.
We have the following result.
{\prop
For $d,k\leq n$, equality
${\mathcal K}_{res}(n,k,d) = \emptyset$ implies
 ${\mathcal K}(n,k,d) = \emptyset$.
}

\pf Proposition \ref{formula2} implies the following for the weight enumerators $w_C$ and $w_{C^{\perp}}$ of $C$ and  $C^{\perp}$:
\begin{equation}
\begin{array}{l}
  w_{C^{\perp}}(x,y)
  = J(C, C^{\perp})(x, y, 0, 0)\\
  = \frac{1}{2^n} J(C,C^{\bot})(x+y,x+y,x-y,x-y)\\
  = \frac{1}{2^n} 2^{n-k} \sum_{u\in C} (x+y)^{n-w(u)} (x-y)^{w(u)}\\
  = \frac{1}{2^k} w_C(x+y,x-y).
\end{array}
\end{equation}
Formula \eqref{RelationA_Bi} then follows from the definition of the Krawtchuk polynomials.
Other formulas are easy. \qed

\subsection{Computational issues}\label{NumResult}

Suppose that for a triple $(n,k,d)$ we find that the polytope ${\mathcal K}(n,k,d)$ is empty.
Then this 
will imply that the triple $(n,k,d)$ is not feasible.
This computation is done by using linear programming iteratively for finding which inequalities
imply equalities when combined. It is time-intensive, which explains why the computation is limited to
$n\leq 16$. We used \cite{cdd} with exact arithmetic as a library in a {\tt C++} program.

If the dimension of ${\mathcal K}(n,k,d)$ is low, then we can compute everything about it; for example,  the facets, vertices and integral points.
The computation is done using \cite{cdd} and direct enumeration of integral points in
an hypercube. Curiously, we found that \cite{4ti2} takes more time than this direct approach.

\subsection{Numerical results}

If $d < d'$, then ${\mathcal K}(n,k,d')\subset {\mathcal K}(n,k,d)$. So, for a given pair
$(n,k)$ there exist a maximum value $d_{max}(n,k)$, such that LCD code with parameters $(n,k,d)$
and $d > d_{max}(n,k)$ are not feasible and ones with $d\leq d_{max}(n,k)$ are not excluded by
the joint LCD enumerator linear programming bound.
The values of $d_{max}$ are given in Table \ref{TheTableD} for $n\leq 16$.
Then for a triple $(n,d)$ we compute $k_{max}(n,d)$ which is the maximum value $k_{max}(n,d)$,
such that for $k > k_{max}(n,d)$, LCD codes of parameters $(n,k,d)$ are not feasible.
The value of $k_{max}$ are given in Table \ref{TheTableK}.
Based on our results, we can state the following conjecture:
{\conj
If ${\mathcal K}(n,k,d)\not=\emptyset$, then ${\mathcal K}(n,k',d)\not=\emptyset$ for
 all $k'<k$.
}

As one can expect, the computation are harder to do when $d$ is small.
The cases $d=1$ and $d=2$ are especially difficult computationally, but, fortunately, for those
cases we can use the trivial codes and the parity check codes.
The {\em parity check code} in dimension $n$ is defined as the  span $C$ of the
vectors $v_1=(1,1,0,\dots, 0)$, $v_2=(0,1,1,0,\dots)$, \dots, $v_{n-1}=(0, \dots, 0,1,1)$.
The 
dual code is $C^{\perp}=(1, \dots, 1)$. The intersection $C\cap C^{\perp}$ is empty if
and only if $n$ is odd. The parameters of the code are $k=n-1$, $d=2$. By using this code
and easy arguments, we can resolve the feasibility of all cases with $d=2$.

\begin{table}
\caption{Some parameters $(n,k,d)$, for which the polytope ${\mathcal K}(n,k,d)$ has low dimension. We give its dimemsion $dim$, number of facets $F$, number of vertices $V$ and number of integral points $P$}\label{LowDimCases}
\begin{center}
\begin{tabular}{||c|c||c|c||}
\hline
$(n,k,d)$ & $(dim, F,V,P)$ & $(n,k,d)$ & $(dim, F,V,P)$\\
\hline
(3,2,2) & (0,1,1,1) & (5,4,2) & (0,1,1,1)\\
(7,2,4) & (0,1,1,1) & (7,3,3) & (0,1,1,1)\\
(7,6,2) & (0,1,1,1) & (8,2,5) & (0,1,1,1)\\
(8,4,3) & (0,1,1,1) & (9,2,6) & (0,1,1,1)\\
(9,4,4) & (0,1,1,1) & (9,8,2) & (0,1,1,1)\\
(10,3,5) & (0,1,1,1) & (11,6,4) & (0,1,1,1)\\
(11,7,3) & (0,1,1,1) & (11,10,2) & (0,1,1,1)\\
(13,2,8) & (0,1,1,1) & (13,12,2) & (0,1,1,1)\\
(14,2,9) & (0,1,1,1) & (15,2,10) & (0,1,1,1)\\
(15,14,2) & (0,1,1,1) & (16,8,5) & (0,1,1,1)\\
(17,3,9) & (0,1,1,1) & (17,8,6) & (0,1,1,1)\\
(4,2,2) & (1,2,2,2) & (5,3,2) & (1,2,2,1)\\
(6,2,3) & (1,2,2,2) & (9,2,5) & (1,2,2,2)\\
(9,5,3) & (1,2,2,2) & (10,2,6) & (1,2,2,2)\\
(10,6,3) & (1,2,2,2) & (12,2,7) & (1,2,2,2)\\
(16,2,10) & (1,2,2,2) & (5,2,2) & (2,4,4,4)\\
(8,2,4) & (3,4,4,3) & (11,2,6) & (3,5,6,5)\\
(12,3,6) & (3,5,5,1) & (17,2,10) & (3,5,6,5)\\
(13,4,6) & (4,7,10,2) & (6,4,2) & (5,12,18,4)\\
(7,2,3) & (5,8,10,3) & (10,2,5) & (5,12,14,7)\\
(13,2,7) & (5,8,10,3) & (14,5,6) & (6,17,121,2)\\
(6,3,2) & (7,19,98,5) & (9,3,4) & (8,16,89,1)\\
(14,3,7) & (8,9,9,1) &  & \\
\hline
\end{tabular}
\end{center}

\end{table}

Also, in Table \ref{LowDimCases} we present several low-dimensional cases.
This Table justifies
 the following conjecture:
{\conj
If the polytope ${\mathcal K}(n,k,d)$ is not empty, then it contains an integral point.
}

In other words, we could not rule out feasibility of $(n,k,d)$-tuples by computing the integral
points. 
But if the polytope contains just one integral point, the possibly feasible
codes have their $M(i,j,k,l)$ values determined and this could be used for further studies.

\begin{table*}[tp]
\caption{Valuee of $d_{max}(n,k)$ for $k\leq n\leq 16$. In parenthesis, best upper bound according to formulation of \cite{DKOSS} if diferent}\label{TheTableD}
\begin{center}
\begin{tabular}{||c||c|c|c|c|c|c|c|c|c|c|c|c|c|c|c|c|}
\hline
n/k & 1 & 2 & 3 & 4 & 5 & 6 & 7 & 8 & 9 & 10 & 11 & 12 & 13 & 14 & 15 & 16\\
\hline
\hline
1 & 1 &  &  &  &  &  &  &  &  &  &  &  &  &  &  & \\
2 & 1 & 1 &  &  &  &  &  &  &  &  &  &  &  &  &  & \\
3 & 3 & 2 & 1 &  &  &  &  &  &  &  &  &  &  &  &  & \\
4 & 3 & 2 & 1 & 1 &  &  &  &  &  &  &  &  &  &  &  & \\
5 & 5 & 2(3) & 2 & 2 & 1 &  &  &  &  &  &  &  &  &  &  & \\
6 & 5 & 3(4) & 2(3) & 2 & 1 & 1 &  &  &  &  &  &  &  &  &  & \\
7 & 7 & 4 & 3(4) & 2(3) & 2 & 2 & 1 &  &  &  &  &  &  &  &  & \\
8 & 7 & 5 & 3(4) & 3 & 2 & 2 & 1 & 1 &  &  &  &  &  &  &  & \\
9 & 9 & 6 & 4 & 4 & 3 & 2 & 2 & 2 & 1 &  &  &  &  &  &  & \\
10 & 9 & 6 & 5 & 4 & 3(4) & 3 & 2 & 2 & 1 & 1 &  &  &  &  &  & \\
11 & 11 & 6(7) & 5(6) & 4(5) & 4 & 4 & 3 & 2 & 2 & 2 & 1 &  &  &  &  & \\
12 & 11 & 7(8) & 6 & 5(6) & 4(5) & 4 & 3(4) & 2(3) & 2 & 2 & 1 & 1 &  &  &  & \\
13 & 13 & 8 & 6(7) & 6 & 5(6) & 4(5) & 4 & 3(4) & 2(3) & 2 & 2 & 2 & 1 &  &  & \\
14 & 13 & 9 & 7(8) & 6(7) & 6 & 5(6) & 4(5) & 4 & 3(4) & 2(3) & 2 & 2 & 1 & 1 &  & \\
15 & 15 & 10 & 7(8) & 6(8) & 6(7) & 6 & 5(6) & 4(5) & 4 & 3(4) & 2(3) & 2 & 2 & 2 & 1 & \\
16 & 15 & 10 & 8 & 7(8) & 6(7) & 6 & 6 & 5 & 4 & 4 & 3 & 2 & 2 & 2 & 1 & 1\\
\hline
\end{tabular}
\end{center}

\end{table*}

\begin{table*}[tp]
\caption{Values of $k_{max}(n,k)$ for $d\leq n\leq 16$. In parenthesis, best upper bound according to formulation of \cite{DKOSS} if different}\label{TheTableK}
\begin{center}
\begin{tabular}{||c||c|c|c|c|c|c|c|c|c|c|c|c|c|c|c|c|}
\hline
n/d & 1 & 2 & 3 & 4 & 5 & 6 & 7 & 8 & 9 & 10 & 11 & 12 & 13 & 14 & 15 & 16\\
\hline
\hline
1 & 1 &  &  &  &  &  &  &  &  &  &  &  &  &  &  & \\
2 & 2 & 0 &  &  &  &  &  &  &  &  &  &  &  &  &  & \\
3 & 3 & 2 & 1 &  &  &  &  &  &  &  &  &  &  &  &  & \\
4 & 4 & 2 & 1 & 0 &  &  &  &  &  &  &  &  &  &  &  & \\
5 & 5 & 4 & 1(2) & 1 & 1 &  &  &  &  &  &  &  &  &  &  & \\
6 & 6 & 4 & 2(3) & 1(2) & 1 & 0 &  &  &  &  &  &  &  &  &  & \\
7 & 7 & 6 & 3(4) & 2(3) & 1 & 1 & 1 &  &  &  &  &  &  &  &  & \\
8 & 8 & 6 & 4 & 2(3) & 2 & 1 & 1 & 0 &  &  &  &  &  &  &  & \\
9 & 9 & 8 & 5 & 4 & 2 & 2 & 1 & 1 & 1 &  &  &  &  &  &  & \\
10 & 10 & 8 & 6 & 4(5) & 3 & 2 & 1 & 1 & 1 & 0 &  &  &  &  &  & \\
11 & 11 & 10 & 7 & 6 & 3(4) & 2(3) & 1(2) & 1 & 1 & 1 & 1 &  &  &  &  & \\
12 & 12 & 10 & 7(8) & 6(7) & 4(5) & 3(4) & 2 & 1(2) & 1 & 1 & 1 & 0 &  &  &  & \\
13 & 13 & 12 & 8(9) & 7(8) & 5(6) & 4(5) & 2(3) & 2 & 1 & 1 & 1 & 1 & 1 &  &  & \\
14 & 14 & 12 & 9(10) & 8(9) & 6(7) & 5(6) & 3(4) & 2(3) & 2 & 1 & 1 & 1 & 1 & 0 &  & \\
15 & 15 & 14 & 10(11) & 9(10) & 7(8) & 6(7) & 3(5) & 2(4) & 2 & 2 & 1 & 1 & 1 & 1 & 1 & \\
16 & 16 & 14 & 11 & 10 & 8 & 7 & 4(5) & 3(4) & 2 & 2 & 1 & 1 & 1 & 1 & 1 & 0\\
\hline
\end{tabular}
\end{center}

\end{table*}

\section{Conclusion}
In this work we have derived a linear programming approach to bounding the parameters of LCD codes.
The linear program considered contains many more variables than the approach in \cite{DKOSS}.
This results in improvements of the upper bound on the minimum distance of LCD codes of given
length and dimension, as evidenced by numerical values.
Bounds on the largest dimension of LCD codes of given length and distance,
which were not considered in \cite{DKOSS} are also given.
Further improvement of the linear programming solvers could allow us to solve larger problems.
It would be interesting to derive a semidefinite approach to these bounds, in the vein
of \cite{SchrijverNewCodeUpperBound} to look for further improvements.

\bibliographystyle{IEEEtran}
\bibliography{IEEEabrv,Ref}

\end{document}